\documentclass{elsart}
\usepackage{amsmath}
\usepackage{amsfonts}
\usepackage{amssymb}
\usepackage[numbers]{natbib}
\usepackage{epic,eepic}
\newcommand{\openbox}{\leavevmode
  \hbox to.77778em{%
  \hfil\vrule
  \vbox to.675em{\hrule width.6em\vfil\hrule}%
  \vrule\hfil}}
\newcommand{\qedsymbol}{\openbox}
\DeclareRobustCommand{\qed}{%
  \ifmmode 
  \else \leavevmode\unskip\penalty9999 \hbox{}\nobreak\hfill
  \fi
  \quad\hbox{\qedsymbol}}
\newenvironment{proof}{\par
  \normalfont\trivlist\item[\hskip\labelsep\itshape Proof.]\ignorespaces
}{\qed\endtrivlist}

\newtheorem{definition}[thm]{Definition}

\makeatletter
\def\@journal{\parbox{0.7\hsize}{\noindent%
}
}
\def\@company{}
\def\ps@copyright{\let\@mkboth\@gobbletwo
  \def\@oddhead{}%
  \let\@evenhead\@oddhead
  \def\@oddfoot{\small\slshape
    \def\@tempa{0}
    \ifx\@volume\@tempa
      \@journal\hfil\@date\/%
    \else
      \@jou@vol@pag\hfil\hbox{}\fi}%
  \let\@evenfoot\@oddfoot
}

\makeatletter

\begin{document}
\begin{frontmatter}
\title{Discrete Nodal Domain Theorems}

\author[KC]{E.\ Brian Davies},
\author[WU]{Josef Leydold}, and
\author[TBI,SFI]{Peter F.\ Stadler}

\address[KC]{Department of Mathematics, King's College,\\
  Strand, London WC2R 2LS, UK\\
  Phone: **44-(0)20 7848 2698 \quad Fax: **44-(0)20 7848 2017 \\
  E-Mail:\quad \texttt{E.Brian.Davies@kcl.ac.uk}\\
  URL: \texttt{http://www.mth.kcl.ac.uk/staff/eb\_davies.html}}

\address[WU]{Dept.\ for Applied Statistics and Data Processing,
  University of Economics and Business Administration,
  Augasse 2-6, A-1090 Wien, Austria\\
  Phone: **43 1 31336-4695 \qquad Fax: **43 1 31336-738\\
  E-Mail:\quad \texttt{Josef.Leydold@statistik.wu-wien.ac.at}\\
  URL: \texttt{http://statistik.wu-wien.ac.at/staff/leydold}}

\address[TBI]{Institute for Theoretical Chemistry and
  Molecular Structural Biology\\
  University of Vienna, W{\"a}hringerstrasse 17, A-1090 Vienna,
  Austria\\
  Phone: **43 1 4277 52737 \qquad Fax: **43 1 4277 52793\\
  E-Mail:\quad \texttt{studla@tbi.univie.ac.at}\\
  URL: \texttt{http://www.tbi.univie.ac.at/\~{ }studla}}

\address[SFI]{The Santa Fe Institute, 1399 Hyde Park Rd,
             Santa Fe NM 87501, USA\\
             E-Mail:\quad\texttt{stadler@santafe.edu}
             }

\begin{abstract}
We give a detailed proof for two discrete analogues of Courant's Nodal
Domain Theorem.

\end{abstract}

\end{frontmatter}

\def\reals{{\mathbb{R}}}
\def\isdef{\,:=\,}
\def\R{\mathbb{R}}

\section{Introduction}

Courant's famous Nodal Domain Theorem for elliptic operators on Riemannian
manifolds (see e.g.\ \cite{Chavel:84}) states\\ {\em If $f_k$ is an
eigenfunction belonging to the $k$-th eigenvalue (written in increasing
order and counting multiplicities) of an elliptic operator, then $f_k$ has
at most $k$ nodal domains.}

When considering the analogous problem for graphs, M.\ Fiedler
\cite{Fiedler:73,Fiedler:75a} noticed that the second Laplacian eigenvalue
is closely related to connectivity properties of the graph, and showed that
$f_2$ always has exactly two nodal domains. It is interesting to note that
his approach can be extended to show that $f_k$ has no more than $2(k-1)$
nodal domains, $k\ge2$ \cite{Powers:88}.  Various discrete versions of the
Nodal Domain theorem have been discussed in the literature
\cite{deVerdiere:93a,Friedman:93a,vanderHolst:96,Duval;Reiner:1999a},
however sometimes with ambiguous statements and incomplete or flawed
proofs. The purpose of this contribution is not to establish new theorems
but to summarize the published results in a single theorem and to present a
detailed, elementary proof.

\section{Preliminaries}

Consider a simple, undirected, loop-free graph $\Gamma$ with finite vertex
set $V$ and edge set $E$.  We write $N\isdef\vert V\vert$ and $x\sim y$ if
$\{x,y\}\in E$. We introduce a weight function $b$ on the edges of
$\Gamma$, conveniently defined as $b:V\times V\to\R$ such that
$b(x,y)=b(y,x)>0$ if $\{x,y\}\in E$ and $b(x,y)=0$ otherwise, and a
potential $v:V\to\R$.  We will consider the {\em Schr{\"o}dinger operator}
\begin{equation}
\mathcal{H}f(x) \isdef \sum_{y\sim x} b(x,y)\left[f(x)-f(y)\right]
                          + v(x)f(x)\,.
\end{equation}
We shall assume that $\Gamma$ is connected throughout this contribution.

The Perron-Frobenius theorem implies that the first eigenvalue $\lambda_1$
of $\mathcal{H}$ is non-degenerate and the corresponding eigenfunction
$f_1$ is positive (or negative) everywhere. Let
\begin{equation}
\lambda_1 < \lambda_2 \le \lambda_3 \le \dots \le \lambda_{k-1} \le
   \lambda_k \le \lambda_{k+1} \le \dots \le \lambda_N
\end{equation}
be the list of eigenvalues of $\mathcal{H}$ arranged in non-decreasing
order and repeated according to multiplicity. Given $k$ let $\overline{k}$
and $\underline{k}$ be the largest and smallest number $h$ for which
$\lambda_h=\lambda_k$, respectively. Let $f_k$ be {\it any} eigenfunction
associated with the eigenvalue $\lambda_k$. Without loss of generality we
may assume that $\{f_i\}$ is a complete orthonormal set of eigenfunctions
satisfying $\mathcal{H}f_i=\lambda_i f_i$.  Since $\mathcal{H}$ is a real
operator, we can take all eigenfunctions to be real.

In the continuous setting one defines the nodal set of a continuous
function $f$ as the preimage $f^{-1}(0)$.  The nodal domains are the
connected components of the complement of $f^{-1}(0)$. In the discrete case
this definition does not make sense since a function $f$ can change sign
without having zeroes. Instead we use the following
\begin{definition}
$D$ is a {\em weak nodal domain} of a function $f:V\to\R$ if it is a
maximal subset of $V$ subject to the two conditions
\begin{itemize}
\item[(i)] $D$ is connected (as an induced subgraph of $\Gamma$);
\item[(ii)] if $x,y\in D$ then $f(x)f(y)\ge0$.
\end{itemize}
$D$ is a {\em strong nodal domain} if (ii) is replaced by
\begin{itemize}
\item[(ii')] if $x,y\in D$ then $f(x)f(y)>0$.
\end{itemize}
\end{definition}

In this contribution we are only interested in nodal domains of
eigenfunctions $f_k$ of the Schr{\"o}dinger operator $\mathcal{H}$.
In the following, the term ``nodal domain'' will always refer to
this case.

The following properties of weak nodal domains are elementary:
\begin{itemize}
\item[(a)] Every point $x\in V$ lies in some weak nodal domain $D$.
\item[(b)] If $D$ is a weak nodal domain then it contains at least one
  point $x\in V$ with $f_k(x)\ne0$ and $f_k$ has the same sign
  on all non-zero points in $D$. Thus each weak nodal domain can be
  called either ``positive'' or ``negative''.
\item[(c)] If two weak nodal domains $D$ and $D'$ have non-empty intersection
  then $f_k\vert_{D\cap D'}=0$ and $D,D'$ have opposite sign.
\end{itemize}
Note that (a) need not hold for strong nodal domains, and (c) is replaced
by: The intersection of two distinct strong nodal domains is empty.

\section{Weak and Strong Nodal Domain Theorem}

The main result of this contribution is
\begin{thm}[Nodal Domain Theorem]
The eigenfunction $f_k$ has at most $\underline{k}$ weak nodal domains and
at most $\overline{k}$ strong nodal domains.
\end{thm}

\begin{proof}
The proof of the Nodal Domain Theorem is based upon deriving a
contradiction from\\
{\bf Hypothesis W:} {\em $f_k$ has $k'>\underline{k}$ weak nodal
domains}, and\\
{\bf Hypothesis S:} {\em $f_k$ has $k'>\overline{k}$ strong nodal
domains},\\
respectively.
\smallskip

We call the domains $D_1,D_2,\dots,D_{k'}$ and define
\begin{equation}
g_i(x) \isdef \left\{
  \begin{matrix} f_k(x) &\text{if }x\in D_i \\
                 0      &\text{otherwise}
  \end{matrix} \right.
\end{equation}
for $1\le i\le k'$. None of the functions $g_i$ is identically zero. Since
they have disjoint supports their linear span has dimension $k'$. It
follows that there exist constants $\alpha_i\in\R$ such that
\begin{equation}
  g\isdef \sum_{i=1}^{k'} \alpha_i g_i
\end{equation}
is non-zero and satisfies $\langle g,f_j\rangle=0$ for $i\le j<k'$. Without
loss of generality we can assume $\langle g,g\rangle=1$, where
$\langle\,.\,,\,.\,\rangle$ denotes the standard scalar product on
$\R^N$. Therefore we have
\begin{equation}
\label{eq:Hgg}
\langle\mathcal{H}g,g\rangle\ge\lambda_{k'}\,.
\end{equation}
Under hypothesis W we know that
\begin{equation}
\label{eq:HggW}
\lambda_{k'}\ge\lambda_k.
\end{equation}
Under hypothesis S we have
\begin{equation}
\label{eq:HggS}
\lambda_{k'}> \lambda_k
\end{equation}
since the last eigenvalue that is equal to $\lambda_k$ has index
$\overline{k}$.

It will be convenient to introduce
$S\isdef\{x\in V\,\vert\, f_k(x)\ne 0\}$
and to define $\alpha:\, V\to\R$ by
\begin{equation}
  \alpha(x)\isdef \left\{
  \begin{matrix} \alpha_i &\mbox{if }x\in S\cap D_i &\mbox{ for some }i \\
           0        &\mbox{otherwise}          & \end{matrix}
                   \right.
\end{equation}
so that $g(x)=\alpha(x)f_k(x)$ for all $x\in V$.

\par\noindent{\bf Lemma 1.} Assuming hypotheses W or S,
we have $\langle\mathcal{H}g,g\rangle\le\lambda_k$.
\par\noindent{\bf Proof.} We have
\begin{equation}
\begin{split}
g(x)\mathcal{H}g(x)
 &= g(x)\sum_{y\sim x} b(x,y)\left[g(x)-g(y)\right] + g^2(x) v(x) \\
 &= \alpha(x)f_k(x)\sum_{y\sim x} b(x,y)\left[
    \alpha(x)f_k(x)-\alpha(y)f_k(y)\right] +\alpha^2(x)f_k^2(x)v(x) \\
 &= \alpha^2(x)f_k(x)\sum_{y\sim x} b(x,y)\left[f_k(x)-f_k(y)\right]
    +\alpha^2(x)f_k^2(x)v(x) \\
 &  \qquad + \alpha(x)f_k(x)\sum_{y\sim x} b(x,y)
         \left[\alpha(x)-\alpha(y)\right]f_k(y)      \\
 &= \alpha^2(x)f_k(x)\mathcal{H}f_k(x) + \mbox{Rem}(x) =
    \alpha^2(x)\lambda_k f_k^2(x)   + \mbox{Rem}(x)\\
 &= \lambda_k g^2(x) + \mbox{Rem}(x)
\end{split}
\end{equation}
Summing over the vertex set yields
\begin{equation}
  \langle \mathcal{H}g,g\rangle = \lambda_k+\mbox{Rem}
\end{equation}
where
\begin{equation}
\begin{split}
\mbox{Rem} &=
   \sum_{x\in V}\sum_{y\sim x}
      b(x,y)\alpha(x)\left[\alpha(x)-\alpha(y)\right]f_k(x)f_k(y) \\
  &= {1\over 2}\sum_{x,y\in V}
      b(x,y)\left[\alpha(x)-\alpha(y)\right]^2f_k(x)f_k(y) \\
\end{split}
\end{equation}
by symmetrizing. A term of the remainder $\mbox{Rem}$ vanishes
if $f_k(x)=0$ or $f_k(y)=0$. If $f_k(x)f_k(y)>0$ and $x\sim y$,
i.e.\ $b(x,y)>0$, then $x$ and $y$ lie in the same nodal domain
and thus $\alpha(x)=\alpha(y)$, and the corresponding contribution
to $\mbox{Rem}$ vanishes as well. The only remaining terms are
those for which $f_k(x)f_k(y)<0$ and $x\sim y$. So we see that
$\mbox{Rem}\le 0$.\\
Thus we have $\langle \mathcal{H}g,g\rangle\le \lambda_k\langle g,g\rangle
=\lambda_k$.
\hfill$\triangle$

Under hypothesis S, eqns.(\ref{eq:Hgg}), (\ref{eq:HggS}), and Lemma~1 lead
to the desired contradiction, proving the second part of the theorem.
\smallskip

Under hypothesis W, eqns.(\ref{eq:Hgg}), (\ref{eq:HggW}), and Lemma~1
imply $\langle g,\mathcal{H}g\rangle=\lambda_k$. Since $g$ is by
construction orthogonal to all eigenvectors $f_j$, $j<\underline{k}<k'$, 
a simple variational argument implies
\begin{equation}
\label{eq:eigenf}
\mathcal{H}g=\lambda_k g\,.
\end{equation}
For the second step of the proof of the Weak Nodal Domain Theorem we
exploit the fact that the remainder $\mbox{Rem}=0$ as a consequence of
equ.(\ref{eq:eigenf}). We proceed with a unique continuation result for the
function $\alpha$.

\par\noindent{\bf Lemma~2.} If hypothesis W holds, $\alpha_i\ne0$,
$x\in D_i$, $y\in D_j\setminus D_i$, and $\{x,y\}\in E$ then
$\alpha_j=\alpha_i$.
\par\noindent{\bf Proof.}
If $x\in D_i$, $y\in D_j\setminus D_i$, $x\sim y$, and $f_k(x)\ne 0$ then
$f_k(y)\ne 0$ (otherwise $y\in D_i\cap D_j$), and hence
$f_k(x)f_k(y)<0$. From $\mbox{Rem}=0$,  $f_k(x)f_k(y)<0$,
and $x\sim y$ we conclude that $\alpha(x)=\alpha(y)$ and hence
$\alpha_i=\alpha(x)=\alpha(y)=\alpha_j$.\\
Now assume that $f_k(x)=0$. Define $h\isdef f_k-(1/\alpha_i)g$.
Then
\begin{equation}
\mathcal{H} h = \lambda_k h\qquad\mbox{and}\qquad
   h\vert_{D_i} = 0\,.
\end{equation}
We have
\begin{equation}
\begin{split}
 0 = \lambda_kh(x) = \mathcal{H}h(x) &=
    \sum_{y\sim x} b(x,y)\left[ h(x)-h(y)\right] + v(x)h(x) \\
   &= -\sum_{y\in B} b(x,y)h(y)
\end{split}
\end{equation}
where $B\isdef \{y\in V\,\vert\, y\sim x\mbox{ and }y\notin D_i\}$.
Note that $B\ne\emptyset$ by the assumptions of the lemma.
Suppose for definiteness that $D_i$ is a positive nodal domain.
Then $y\in B$ satisfies $f_k(y)<0$ since otherwise one would have
to adjoin $y$ to $D_i$. Thus $B\cup\{x\}$ is a connected set on which
$f_k\le0$. Therefore it is contained in the single (negative) nodal
domain $D_j$. Therefore
\begin{equation}
0 = -\sum_{y\in B} b(x,y)h(y) =
 - \left(1-{\alpha_j\over\alpha_i}\right) \sum_{y\in B} b(x,y)f_k(y).
\end{equation}
The terms in the sum are all negative, thus $\alpha_i=\alpha_j$.\\
The same argument of course works when $D_i$ is a negative nodal domain.
\hfill$\triangle$

We say that $D_i$ {\em is adjacent to} $D_j$ if there are $x\in D_i$ and
$y\in D_j\setminus D_i$, $x\sim y$. Note that adjacent nodal domains
must have opposite signs. Now consider a collection $\{D_1,\dots,D_l\}$
of nodal domains such that $\bigcup_i D_i\ne V$. Then there exists a
nodal domain $D_j\ne D_i$, $i=1,\dots,l$, that is adjacent to some $D_i$,
$i=1,\dots,l$; otherwise $\Gamma$ would not be connected.

Now we are in the position to prove the first part of the theorem.
We assume hypothesis W and thus the conclusions of lemma~1 and lemma~2.
Since $g\ne 0$ there exists an index $i$ for which $\alpha_i\ne0$.
If $D_j$ is a nodal domain adjacent to $D_i$ then lemma~2 implies
$\alpha_j=\alpha_i$. Since the graph $\Gamma$ is connected by
assumption, we conclude in a finite number of steps that
$\alpha_j=\alpha_i$ for all $j$. Hence $g=\alpha_i f_k$.
This, however, contradicts the fact that $\langle g,f_k\rangle=0$.
\end{proof}

\section{Two Counter-Examples}

Neither the Weak nor the Strong Nodal Domain theorem can be strengthened
without additional assumptions. If $\Gamma$ is a path with $N$ vertices,
then $f_k$ has always $k$ weak nodal domains.  An example where $f_k$ has
more than $k$ strong nodal domains is e.g.\ given by Friedman
\cite{Friedman:93a}: a star on $n$ nodes, i.e., a graph which is a tree
with exactly one interior vertex, has a second eigenfunction with $n-1$
strong nodal domains. For example, the star with $5$ nodes has
$\lambda_2=\lambda_3=\lambda_4=1$ and an eigenvector $f_2=(0,1,1,-1,-1)$,
where the first coordinate refers to the interior vertex. Since $f_2$
vanishes at the interior vertex each of the $n-1$ leafs is a strong nodal
domain.  These eigenvectors of the stars may also serve as a counterexample
to Theorem~6 and Corollary~7 of \cite{Duval;Reiner:1999a}.

Theorems~2.4 of \cite{Friedman:93a} and 4.4 of \cite{vanderHolst:96} can
be rephrased as follows: {\em If $f_k$ has more than $k$ strong nodal
domains, then there is no pair of vertices such that $f_k(x)>0$,
$f_k(y)<0$ and $x\sim y$, i.e., there is no edge that joins any two
strong nodal domains.} This statement is incorrect, as the following
example shows:

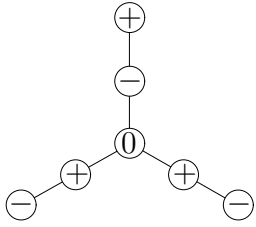
\begin{center}
\setlength{\unitlength}{0.0004in}
\begin{picture}(3038,2664)(0,-10)
\path(1523,2536)(1523,908)\path(1530,901)(90,90)\path(1538,901)(2955,90)
\put(810,490){\whiten\ellipse{400}{400}}
\put(670,400){$+$}
\put(90,90){\whiten\ellipse{400}{400}}
\put(-70,0){$-$}
\put(2235,490){\whiten\ellipse{400}{400}}
\put(2090,400){$+$}
\put(2955,90){\whiten\ellipse{400}{400}}
\put(2805,0){$-$}
\put(1523,1726){\whiten\ellipse{400}{400}}
\put(1370,1620){$-$}
\put(1523,2566){\whiten\ellipse{400}{400}}
\put(1370,2460){$+$}
\put(1523,908){\whiten\ellipse{400}{400}}
\put(1410,780){$0$}
\end{picture}
\end{center}

This tree has eigenvalues $\lambda_5=\lambda_6=(3+\sqrt{5})/2$ and
a corresponding eigenvector
\begin{equation}
f_5 = (2, -1-\sqrt{5}, 0, (1+\sqrt{5})/2, (1+\sqrt{5})/2, -1, -1)
\end{equation}
from top to bottom. There are $5$ weak and $6$ strong nodal domains.
Nevertheless, there are edges connecting strictly positive with strictly
negative vertices.

\section*{Acknowledgements}

This work was supported in part by the Austrian {\it Fonds zur
F{\"o}rderung der Wissenschaftlichen Forschung}, Proj.\ No.\ P14094-MAT.


\end{document}